\documentclass[10pt]{article} 
\usepackage{amsfonts,amssymb,latexsym,amsmath,oldgerm,amscd}
\newcommand{\mf}{\mathfrak}

\newcommand{\lra}{\longrightarrow}

\newcommand{\ol}{\overline}

\newcommand{\ra}{\rightarrow}
\newcommand{\Ra}{\Rightarrow}

\newcommand{\mbb}{\mathbb}
\newcommand{\tn}{\textnormal}

\newtheorem{de}{Definition}[section]
\newtheorem{re}[de]{Remark}
\newtheorem{pr}[de]{Proposition} 
\newtheorem{tr}[de]{Theorem}
\newtheorem{lm}[de]{Lemma} 

\newtheorem{nt}[de]{Notation} 

\newtheorem{co}[de]{Corollary}

\def\vp{\rm \vspace{0.2cm}}

\def\hb{\hfill$\Box$}
\def\aut{\rm Aut} 
\def\sp{\rm Spec} 
\def\tran{\rm Trans}
\def\fl{\rm ETrans}
\def\Um{\rm Um}
\def\GL{\rm GL}
\def\SL{\rm SL}
\def\EO{\rm EO}
\def\SO{\rm SO}
\def\E{\rm E}
\def\T{\rm T}
\def\ET{\rm ET}
\def\G{\rm G}
\def\Sp{\rm Sp}

\def\k{\rm K_1}
\def\K{\rm K}
\def\T{\rm T}
\def\O{\rm O}
\def\ESp{\rm ESp}
\def\EO{\rm EO}
\def\es{\rm S}
\def\sdim{\rm sdim}
\begin{document}
\title{Injective Stability for ${\k}$ of Classical Modules}
\author{ Rabeya Basu \& Ravi A. Rao\\
{\small Tata Institute of Fundamental Research, Homi Bhabha Road, 
Mumbai-400005, India}
\footnote{The first author was partially supported by the T.I.F.R. 
Endowment Fund.}
\footnote{Correspondence author: Ravi A. Rao; 
{\it email: ravi@math.tifr.res.in}}}
\date{}
\maketitle
\begin{center}{\it 2000 Mathematics Subject Classification:
{13C10, 13H05, 15A63, 19B10, 19B14}}
\end{center}
\begin{center}{\it Key words: regular ring, affine algebra, 
projective modules, ${\k}$, ${\k}{\Sp}$}
\end{center} {\small ~~~~Abstract: In \cite{RV}, the second author and
W. van der Kallen showed that  the injective stabilization bound for
${\k}$ of general linear group is $d+1$ over a regular affine algebra
over a perfect $C_1$-field, where $d$ is the krull dimension of the
base ring and it is  finite and at least 2. In this article we prove
that the injective stabilization bound  for ${\k}$ of the symplectic group
is $d+1$ over a geometrically regular ring containing a field, where
$d$ is the stable dimension of the base ring and it is finite and at
least 2. Then using the Local-Global Principle for the transvection
subgroup of the automorphism group of projective and symplectic
modules we show that the injective stabilization bound is $d+1$ for
${\k}$ of projective and symplectic modules of global rank at least 1
and local rank at least 3 respectively in each of the two cases above.}

\section{\large Introduction} 

~~~~In this article we discuss the injective stabilization 
for the ${\k}$ group of projective and symplectic  modules.

In the early 1960's Bass-Milnor-Serre began the study of the stabilization
for the linear group ${\GL}_n(R)/{\E}_n(R)$ for $n\ge 3$, where $R$ is a 
commutative ring with identity. In \cite{B3}, they showed that 
${\k}(R)={\GL}_{d+3}(R)/{\E}_{d+3}(R)$, where $d$ is the dimension of the
maximum spectrum. (They also showed that ${\k}(R) = {\GL}_3(R)/{\E}_3(R)$, when Krull dimension of $R$ is $1$.) 
In \cite{V}, L.N. Vaserstein proved their conjectured 
bound of $(d+2)$ for an associative ring with identity, where $d$ is the 
stable dimension of the ring. After that, in \cite{V2}, he studied the 
orthogonal and the 
unitary ${\k}$-functors, and obtained stabilization theorems for them. 
He showed that the natural map 
$$\begin{cases} \varphi_{n,n+1}: \frac{{\es}(n,R)}{{\E}(n,R)}\lra 
\frac{{\es}(n+1,R)}{{\E}(n+1,R)} & \mbox{ in the linear case } \\
\varphi_{n,n+2}:\frac{{\es}(n,R)}
{{\E}(n,R)}\lra \frac{{\es}(n+2,R)}{{\E}(n+2,R)}  
& \mbox{ otherwise } 
\end{cases}$$
(where ${\es}(n,R)$ is the group of automorphisms of the projective,
symplectic and orthogonal modules of rank $n$ with determinant $1$, and 
${\E}(n,R)$ is the elementary subgroup in the respective cases)
is surjective for $n\ge d+1$ 
in the linear case, for $n\ge d$ in the symplectic case, and for
$n\ge 2d+2$ in the orthogonal case, and is injective for  
$n\ge 2d+4$ in the symplectic and the orthogonal cases. 
Soon after, in \cite{V3}, he
studied stabilization for groups of automorphisms of modules over rings and
modules with quadratic forms over rings with involution, and obtained 
similar stabilization results. 

In \cite{RV}, the second author and W. van der Kallen showed that 
if $A$ is a non-singular affine algebra of dimension $d>1$ 
over a perfect $C_1$-field (Definition 4.1), then the 
natural map 
$$\frac{{\SL}_n(A)}{{\E}_n(A)}\lra \frac{{\SL}_{n+1}(A)}{{\E}_{n+1}(A)}$$ 
is injective for $n\ge d+1$. We generalize this result for the 
automorphism group of finitely generated projective module of 
global rank at least 1 and local rank at least 3. 
(By definition the global 
rank or simply rank of a finitely generated projective $R$-module 
(resp. symplectic or orthogonal $R$-module) 
is the largest integer 
$r$ such that $\overset{r}\oplus R$ (resp. $\overset{r}\perp \mbb{H}(R)$) 
is a direct summand (resp. orthogonal summand) of the module. 
$\mbb{H}(R)$ denotes the hyperbolic plane). More precisely, 
we prove the following: (We assume that $(H1)$ and 
$(H2)$ holds, as stated in \ref{note}). 
\vp \\

{\bf Theorem 1.}
\label{st1} 
{\it Let $A$ be an affine algebra of dimension $d>1$ over a perfect 
$C_1$-field $k$. Assume $(d+1)\,!A=A$.
Let $P$ be a finitely generated projective $A$-module of local rank $n>1$.
If $\gamma\in {\SL}(P)$ is such that 
$\gamma\perp \{1\}\in {\tran}(P\oplus A)$ 
and $n\ge d+1$, then $\gamma$ is isotopic to the identity,
{\it i.e.} there exists an automorphism 
$\alpha(X)\in {\SL}(P[X])$ such that $\alpha(0)=\tn{Id}$ and 
$\alpha(1)=\gamma$. Moreover, if $A$ is non-singular, then 
$\gamma\in {\tran}(P)$. In particular, 
the map $\rho: \frac{{\SL}(P)}{{\tran}(P)}\lra 
\frac{{\SL}(P\oplus A)}{{\tran}(P\oplus A)}$ 
is bijective for $n\ge d+1$.} \vp

{\bf Theorem 2.}
{\it Let $R$ be a commutative ring with identity of stable
dimension $d>1$ and $A$ be an associative $R$-algebra such that $A$ is finite 
as a left $R$-module. Let $(P,\langle \, , \rangle)$ be a symplectic left
$A$-module of even local rank $n\ge \tn{max }(3,d+1)$. 
If $\gamma\in {\Sp}(P)$ is such that 
$\gamma\perp I_2\in {\tran}_{\Sp}(P\perp A^2)$, then $\gamma$ is 
isotopic to the identity, {\it i.e.} there exists an automorphism 
$\alpha(X)\in {\Sp}(P[X])$ such that $\alpha(0)=\tn{Id}$ and 
$\alpha(1)=\gamma$. Moreover, if $A$ is a geometrically regular ring 
containing a field $k$, then 
$\gamma\in {\tran}_{\Sp}(P)$. In particular, the map 
{\small $$\rho_{\Sp}:\frac{{\Sp}(P)}
{{\tran}_{\Sp}(P)}\lra 
\frac{{\Sp}(P\perp A^2)}{{\tran}_{\Sp}(P\perp A^2)}$$} 
is bijective for $n\ge \tn{max} (3,d+1)$. } \vp

However, in a companion article \cite{brj} 
we prove that the injective stabilization bound for ${\k}$ of the orthogonal group 
is not less than $2d+4$, in general, for an affine algebra over a 
perfect $C_1$-field.\vp

\section{\large Preliminaries} 
\begin{de}\tn{Let $R$ be an associative ring with identity.
The following condition was introduced by H. Bass in \cite{B2}:}

\tn{{\bf $(R_m)$} for every $(a_1,\dots, a_{m+1})
\in {\Um}_{m+1}(R)$, there are $\{x_i\}_{(1\le i\le m)}\in R$ such that} 
$(a_1+a_{m+1}x_1)R+\cdots+(a_m+a_{m+1}x_m)R=R$.

\tn{The condition $(R_m)\Ra (R_{m+1})$ for every $m> 0$. Moreover, for any 
$n\ge m+1$ the condition $(R_m)$ implies $(R_n)$ with $x_i=0$ for $i\ge m+1$.}

\tn{By {\bf stable range} for an associative ring $R$  
we mean the least $n$ such that $(R_n)$ holds.}

\tn{Although, it appears that we should have referred
to the above condition as $R$ having ``right'' 
stable range $n$, it has been shown by L.N. Vaserstein 
(\cite{V1}, Theorem 2) that 
``right stable range $n$'' and ``left stable range $n$'' are 
actually equivalent conditions. The integer $n-1$ is called the {\bf stable 
dimension} of $R$ and is denoted by ${\sdim}(R)$.}
\end{de} 
\begin{lm} \label{sd1} $({\it cf. } \cite{B2})$
If $R$ is a commutative noetherian ring with identity of Krull dimension $d$, 
then ${\sdim}(R)\le d$.
\end{lm}

\begin{de} \tn{Let $R$ be an associative ring with identity.
To define other classical modules, we need an involutive antihomomorphism
({\bf involution}, in short) $*:R\ra R$ ({\it i.e.},  
$(x-y)^{*}=x^{*}-y^{*}$, $(xy)^{*}=y^{*}x^{*}$ and $(x^{*})^{*}=x$ for any 
$x,y\in R$. We assume that $1^{*}=1$. 
For any left $R$-module $M$ the involution induces a left module 
structure to the right $R$-module $M^{*}$=Hom$(M,R)$ given by 
$(xf)v=(fv)x^{*}$, where $v\in M$, $x\in R$ and $f\in M^{*}$. In this case if
$M$ is a left $R$-module then $O_M(m)$ has a right $R$-module structure. 
But any right $R$ module can be viewed as a left $R$-module via the 
convention $ma=a^{*}m$ for $m\in M$ and $a\in R$. }
\end{de} 
{\bf Blanket Assumption:} Let $A$ be an $R$-algebra, 
where $R$ is a commutative ring with identity, such that $A$ is finite as a 
left $R$-module.
Let $A$ possess an involution 
$*:r\mapsto \bar{r}$, for $r\in A$. For a matrix $M=(m_{ij})$ over $A$ we  
define $\ol{M}=(\ol{m}_{ij})^t$. Let 
$\psi_1=\begin{pmatrix} 0 & 1 \\ -1 &0 \end{pmatrix}$, 
$\psi_n=\psi_{n-1}\perp \psi_1$, for $n>1$.
For a column vector $v\in A^n$ we write 
$\widetilde{v}=\bar{v}^t\cdot \psi_n$ in the symplectic case.

We define a form $\langle\, ,\rangle$ as follows: 
$$\langle v,w\rangle =\begin{cases} v^t\cdot w  & \mbox{in the linear case} \\
\widetilde{v}\cdot w & \mbox{in the symplectic case.} \end{cases}$$
(Viewing $M$ as a right $A$-module we can assume the linearity).

Since $R$ is commutative, we can assume that 
the involution ``$*$'' defined on $A$ is trivial over 
$R$. We shall always assume that 2 is invertible in the ring $R$ while 
dealing with the symplectic case. For definitions of 
the automorphism group, the symplectic  module, and its transvection and its 
elementary transvection subgroup, see (\cite{BR2}, \S 2).

\begin{nt} \label{note}
\tn{In the sequel $P$ will denote either a finitely generated 
projective left $A$-module of local rank 
$n$, a symplectic left module of even rank 
$n=2r$ with a fixed form $\langle \, ,\rangle$. And 
$Q$ will denote $P\oplus A$ in the linear case and $P\perp A^2$
in the symplectic case. To denote $(P\oplus A)[X]$ in the linear case and 
$(P\perp A^2)[X]$ in the symplectic case we will use the notation $Q[X]$. 
We assume that the local rank of projective module is at least 3 
when dealing with the linear case and at least 6 when considering the 
symplectic case. For a finitely generated projective 
$A$-module $M$ we use the notation 
${\G}(M)$ to denote the automorphism group of the projective module = ${\aut}(M)$ and 
the group of isometries of the symplectic module = ${\Sp}(M,\langle \, ,\rangle )$. 
Let ${\SL}(M)$ denote the automorphism group of the projective module with determinant 1  
in the case when $A$ is commutative. We use
${\es}(M)$ to denote ${\SL}(M)$ in the linear case and ${\Sp}(M,\langle \, ,\rangle )$ in the symplectic case. Let 
${\T}(M)$  denote the transvection subgroup of the automorphism group of the projective module  ${\tran}(M)$,  
and the transvection subgroup of the automorphism group of the symplectic module = ${\tran}_{\Sp}(M)$. We write 
$\tn{ET}(M)$ to denote the elementary transvection subgroup of the automorphism group of the projectivs module  ${\fl}(M)$, and the 
transvection subgroup of the automorphism group of symplectic module ${\fl}_{\Sp}(M)$. (For details see \cite{BR2}).}
\end{nt}  
We shall assume \vp\\
{\bf (H1) for every maximal ideal $\mf{m}$ of $A$, $Q_{\mf{m}}$ is isomorphic 
to $A^{2n+2}_{\mf{m}}$ with the standard bilinear form $\mbb{H}(A_{\mf{m}}^{n+1})$.} 
\vp\\
{\bf (H2) for every non-nilpotent $s\in A$, if the projective module 
$Q_s$ is free $A_s$-module, then the symplectic module 
$Q_s$ is isomorphic to $A^{2n+2}_s$ with the standard bilinear form  
$\mbb{H}(A_s^{n+1})$. }

\begin{nt} \label{note1}
\tn{When $P=A^n$ (n is even is the non-linear cases), we also use the notation
${\G}(n,A)$, ${\es}(n,A)$ and 
${\E}(n,A)$ for ${\G}(P)$, ${\es}(P)$ and ${\T}(P)$ respectively. 
We denote the usual standard elementary generators of ${\E}(n,A)$ by 
$ge_{ij}(x)$, $x\in A$. $e_i$ will denote the column vector  
$(0,\dots,1,\dots,0)^t$ (1 at the i-th position).}
\end{nt}
\begin{re} \label{endo}
\tn{Let $Q$ be as in \ref{note}. Note that if $\alpha\in \tn{End}(Q)$ then 
$\alpha$ can be considered as a matrix of the form
{\small $\left( \begin{array}{cc} \tn{End}(P) & \tn{Hom}(P,A) \\ 
\tn{Hom}(A,P) & \tn{End}(A) \end{array}\right)$} in the linear case. 
In the non-linear cases one has a similar matrix for $\alpha$ of the form 
{\small $\left( \begin{array}{cc} \tn{End}(P) & \tn{Hom}(P,A\oplus A) \\ 
\tn{Hom}(A\oplus A,P) & \tn{End}(A\oplus A) \end{array}\right)$}}.
\end{re}

\begin{de} \tn{An associative ring $R$ is said to be {\bf semilocal} if 
$R/\tn{rad}(R)$ is artinian semisimple.}
\end{de} 
\begin{lm} \label{HB}
\tn{(H. Bass) ({\it cf.} \cite{B2})} 
Let $A$ be an associative $R$-algebra such that $A$ is finite as a left 
$R$-module and $R$ be a commutative local ring with identity. 
Then $A$ is semilocal. 
\end{lm} 
{\bf Proof.} Since $R$ is local, $R/\tn{rad} (R)$ is a division ring by 
definition. That implies $A/\tn{rad} (A)$ is a finite module over the division 
ring $R/\tn{rad}(R)$ and hence is a finitely generated vector space. Thus 
$A/\tn{rad} (A)$ artinian as $R/\tn{rad}(R)$ module and hence  
$A/\tn{rad}(A)$ artinian as $A/\tn{rad}(A)$ module, so 
it is an artinian ring. 
It is known that a right artin ring is semisimple if its radical is trivial. Now $\tn{rad}(A/\tn{rad}(A))=0$, hence it follows that  
$A/\tn{rad}(A)$ is semisimple. 
Hence $A/\tn{rad}(A)$ artinian semisimple. Therefore, $A$ is semilocal by 
definition. \hb \vp 

We recall the well-known Serre's unimodular theorem:

\begin{tr} \tn{(J-P. Serre)} $({\it cf.} \cite{B2})$ Let $R$ be a commutative noetherian ring of 
dimension $d$, and let $P$ be a finitely generated projective $R$-module of 
local rank $\ge d+1$. Then $P$ contains a unimodular element. 
\end{tr}  

While dealing with the symplectic case we implicitly use the following 
well-known fact; which we include for completeness. 

\begin{lm} Let $R$ be a commutative ring with identity and
$(P,\langle\, ,\rangle)$ be a symplectic $R$-module. If $P$ contains a 
unimodular element, 
then $(P,\langle\, ,\rangle)$ contains a hyperbolic plane as a direct summand. 
\end{lm} 
{\bf Proof.} Let $p\in {\Um}(P)$ and let $\varphi:P\cong P^{*}$ be the induced 
isomorphism. Then there exists $\alpha:P\ra R$ such that $\alpha(p)=1$. 
Since $\langle p,p\rangle=0$, it follows that $p\ne \varphi^{-1}(\alpha)$. 
Hence there exists $f\in P$ such that $f\ne p$ and $\varphi(f)=\alpha$. Now if 
$x\in Rp\cap Rf$, then $x=tp=sf$, for some $t,s\in R$. Since 
$\langle x,x\rangle=0$, it follows that $st=0$. Hence $sx=0$. 
This is a contradiction, as $Rp\cong R$. Hence $Rp\cap Rf=0$. 
Also $\langle p,f\rangle=1$; hence $P$ contains $\mbb{H}(R)$. We claim that 
$P$ contains $\mbb{H}(R)$ as a direct summand. Let 
$$Q=\{q\in P \,|\, \langle q,f\rangle=0, \langle q,p\rangle=0\}.$$
Again, let $y\in Q\cap (Rp\oplus Rf)$. Then $y=ap+bf$ for some $a,b\in R$.  
Since $\langle y,p\rangle=\langle y,f\rangle=0$, it will follow that $y=0$. 
Hence $Q\oplus (Rp\oplus Rf)\subseteq P$. Now let $z\in P$ be such that 
$z\ne p$ and $z\ne f$. Let $z'=z-\langle z,p\rangle f+\langle z,f\rangle p$. 
Then one checks that $z'\in Q$. Hence $P\cong Q\oplus \mbb{H}(R)$. 
(Note: $Q$ inherits 
a symplectic structure from $P$ given by the restriction 
$\langle\,,\rangle\,|_Q:Q\times Q\ra R$). Hence the result follows. \hb 

The following theorem is a well known result: 
\begin{tr} \label{tr3}
Let $R$ be an associative ring of stable dimension $d\ge 1$. Then, for
$n\ge d+2$ in the linear case and for $n\ge 2d+4$ in the symplectic and the 
orthogonal cases, ${\E}(n,R)$ acts transitively on ${\Um}_n(R)$. 
In other words, any unimodular row of length $n$ over $R$ is completable 
to an elementary matrix if $n\ge d+2$ in the linear case and 
$n\ge 2d+4$ in the symplectic case. 
\end{tr}
{\bf Proof.} See (\cite{LAM}, Theorem $7.3'$, pg. $93$) for the linear case 
and (\cite{V2}, Theorem 2.7) for the symplectic case. 
(The key to proving it is Lemma \ref{sd1}).\hb

\begin{de} \tn{For $\alpha\in {\rm M}(r,R)$ and $\beta\in {\rm M}(s,R)$ 
we have $\alpha\perp \beta$ denotes its 
embedding $M({r+s},R)$ given by ($r$ and $s$ are even in the non linear cases)}
$$\alpha\perp \beta =
\left(\begin{array}{cc} \alpha & 0 \\ 0 & \beta
\end{array}\right).$$ 
\tn{There is an infinite counterpart: 
Identifying each matrix $\alpha\in {\GL}_n(R)$ 
with the large matrix $(\alpha\perp \{1\})$
gives an embedding of ${\GL}_n(R)$ into ${\GL}_{n+1}(R)$. 
Let ${\GL}(R)=\underset{n=1}
{\overset{\infty}\cup} {\GL}_n(R)$, 
${\SL}(R)=\underset{n=1}{\overset{\infty}\cup} {\SL}_n(R)$, and 
${\E}(R)=\underset{n=1}{\overset{\infty}\cup} {\E}_n(R)$ be the corresponding 
infinite linear groups.}
\end{de}
\begin{de} \tn{The quotient group}
$${\k}(R)=\frac{{\GL}(R)}{[{\GL}(R),{\GL}(R)]}=\frac{{\GL}(R)}{{\E}(R)}$$
\tn{is called the {\bf Whitehead group} of the ring $R$. For }
\tn{$\alpha\in {\GL}_n(R)$ let $[\alpha]$ denote its equivalence class in 
${\k}(R)$. Similarly, one can define the Symplectic Whitehead group 
${\k}{\Sp}(R)$. }
\end{de}

The following theorem is the key result we use to generalize the results known for free 
modules to classical modules. Here we state the result. For details 
see \cite{BR2}. 
\begin{tr} {\bf (Local-Global Principle)} $({\it cf.} \cite{BR2})$\label{lgt}
Let $A$ be an associative $R$-algebra such that $A$ is finite as a left 
$R$-module and $R$ be a commutative ring with identity. 
Let $P$ and $Q$ be as in \ref{note}. Assume that $($H1$)$ holds.
Suppose $\sigma(X)\in {\G}(Q[X])$ with
$\sigma(0)=\tn{Id}$. If 
$$\sigma_{\mf{p}}(X)\in \begin{cases} {\E}(n+1,A_{\mf{p}}[X]) & 
\mbox{in the linear case, }\\
{\E}(2n+2,A_{\mf{p}}[X])& \mbox{in the symplectic case, } 
\end{cases}$$ 
for all $\mf{p}\in {\sp}(R)$, then $\sigma(X)\in {\rm ET}(Q[X])$.
\end{tr} 
\begin{co} \label{tf} 
Let $A$ be an associative $R$-algebra such that $A$ is finite as a left 
$R$-module and $R$ be a commutative ring with identity. Let
$Q$ be as in \ref{note}. Assume that $($H1$)$ holds. 
Then ${\T}(Q)=\tn{ET}(Q)$. 
\end{co}

 \section{\large Stabilization Bounds for ${\k}$ of Classical Modules}

~~~~In this section we prove Theorem 1 and Theorem 2 
stated in Section 1. We will show that the injective stability estimates for 
${\k}(R)$ and ${\k}{\Sp}(R)$, stated by L.N. Vaserstein in 
\cite{V3}, can be improved in the linear and the symplectic cases 
if $R$ is a regular affine algebra over a perfect $C_1$-field. Recall

\begin{de} \tn{ Let $A$ be an affine algebra of dimension 
$d$ over a field $k$ satisfying: For any prime $p\le d$ one of the 
following conditions is satisfied: $(i)$ $p\ne \tn{char } k, c.d._p k\le 1$, 
$(ii)$ $p=\tn{char } p$ and $k$ is perfect. 
In this case we say that $A$ is an 
{\bf affine algebra over a perfect $C_1$-field.}}
\end{de}

Suslin showed that stably free projective modules of top rank $d$ over an affine 
algebra over a field $k$, in which $d!$ was invertble, are free if $k$ is algebraically closed 
in \cite{su1}; and over perfect C$_1$-fields in \cite{su2}. His methods were used to prove their 
cancellative properties in \cite{Bh}; who established the following:

\begin{tr} \label{bhat} 
\tn{(S.M. Bhatwadekar)} $(\cite{Bh}, \tn{Theorem } 4.1)$
Let $A$ be an affine algebra of krull dimension $d>1$ over a perfect $C_1$-field 
$k$. Assume $d\,!A=A$.  Let $P$ be a projective $A$-module 
of local rank $d$. Then for $(p,a)\in \tn{Um}(P\oplus A)$ there exists 
$\tau\in {\aut}(P\oplus A)$ such that $(p,a)\tau=(0,1)$. 
In particular, $P$ is cancellative. 
\end{tr}
\begin{lm} $({\it cf.} \cite{RV})$ \label{rv}
Let $A$ be as in Theorem 
\ref{bhat} and let $P$ be a projective $A$-module of local rank $d+1$, 
where $d$ is the stable dimension of $A$ and $d>1$.
Let $v(X)\in \tn{Um}((P\oplus A)[X])$ with
$v(X)\equiv (0,1)$ modulo $(X^2-X)$. Then there 
exists $\sigma(X)\in {\SL}((P\oplus A)[X])$ with $\sigma(X)\equiv$ {\tn Id} 
modulo $(X^2-X)$ such that $v(X)\sigma (X)=(0,1)$. 
\end{lm} 
{\bf Proof.} Our argument is similar to that in (\cite{RV}, Proposition 3.3).
Let $Y=X^2-X$ and $B=A[Y,Z]/(Z^2-YZ)$. Then $B$ is an affine 
algebra of dimension $d+1$ over the field $k$. Let $v(X)=e_{d+2}^t+Yv'(X)$ 
with $v'(X)\in (P\oplus A)[X]$. Let $u(Z)=e_{d+2}^t+Zv'(X)$ be its lift in 
$B$. Then $u(Z)\in \tn{Um}((P\otimes B)\oplus B)$ as locally it is unimodular. 
So $u(Y)=v(X)$ and $u(0)=(0,1)$. By Proposition \ref{bhat} there exists 
$\beta(Z)\in {\SL}((P\otimes B)[Z]\oplus B[Z])$ with $u(Z)=\beta(Z)e_{d+2}^t$. 
Take $\sigma=\beta(0)^{-1}\beta(Y)$. \hb 
\begin{lm} \label{b1} 
\tn{(\cite{B2}, Chapter IV, Theorem 3.1)} 
Let $A$ be an associative $R$-algebra such that $A$ is finite as a left 
$R$-module where $R$ is a commutative ring with identity.
Let $A$ have stable dimension $d$, 
and let $P$ be a projective left $A$-module of rank $n$, where
$n\ge d+1$. Suppose that $(p,a)\in \tn{Um}(P\oplus A)$. Then 
there exists a homomorphism $f:A\ra P$ such that $p+f(a)\in \tn{Um}(P)$
and there exists $\tau\in {\tran}(P\oplus A)$ such that $\tau:(p,a)\ra (0,1)$. 
\end{lm} 

\begin{lm} \label{b3a} \tn{(L.N. Vaserstein, \cite{V3})}
Let $A$ be a commutative $R$-algebra such that $A$ is 
finite as a $R$-module and $R$ is a commutative ring with identity.
Let $(P,\langle \, , \rangle)$ be a symplectic left $A$-module of even local 
rank $n$, where $n\ge d$. 

Let $(p,b,a)\in \tn{Um}(P\perp A^2)$.
\begin{enumerate} 
\item There exists $\tau\in {\tran}_{\Sp}(P\perp A^2)$ such that 
$\tau:(p,b,a)\ra (0,0,1)$. 
\item If  $I$ is a two sided ideal of $A$ and 
$(p,b,a)\equiv (0,0,1)$ modulo $I(P\perp A^2)$, then  
there exists $\tau\in {\tran}_{\Sp}(P\perp A^2,I)$ such that 
$\tau:(p,b,a)\ra (0,0,1)$. 
\end{enumerate}
\end{lm} 
{\bf Proof.} We prove the result for completeness.  
We follow the line of proof of R.G. Swan, (see \cite{S}, Corollary 9.8) 
(Also see  \cite{Bh}, Theorem 3.2).

By Lemma \ref{b1} there exists $q\in P$ and $t\in A$ such that 
${\O}(p+aq,b+at)= A$; {\it i.e.} $(p+aq,b+at)\in \tn{Um} (P\perp A)$. 
Hence there exists $\gamma\in P^{*}$ and $g\in A$ such that 
$\gamma(p+aq)+g(b+at)=1$. 
Let $\eta=g\Phi(q)$, where $\Phi:P \cong P^{*}$ is the 
induced isomorphism. Then $\eta (p+aq)= -g\langle p,q\rangle$. Hence 
$\delta(p+aq)+g(b+at+\langle p,q\rangle)=1,$ where $\delta=\eta+\gamma$.
Now consider the following automorphisms (elementary transvections) 
of $(P\perp A^2)$: 
$$\theta_{(t,q)}:(p,b,a)\mapsto (p+aq,b+at+\langle p,q\rangle,a),$$ 
$$\tau_{(g,\beta)}:(p,b,a)\mapsto (p-\beta(b),b,a+gb+\delta(p)),$$ where 
$\beta: A\ra P$ with $\beta^{*}=\delta\Phi^{-1}$.   
Let $\tau'_{(g,\beta)}=(1-a)\tau_{(g,\beta)}$,
$\tau=\theta_{(-b_1,-p_1)}\tau_{(g,\beta)}\theta_{(t,q)}$ for 
$b_1=b+ta+\langle p,q\rangle$, and $p_1=p+aq-\beta(b+ta+\langle p,q\rangle)$. 
Then $\tau (p,b,a)^t=(0,0,1)^t$; as required.

Next assume $(p,b,a)\equiv (0,0,1)$ modulo 
$I(P\perp A^2)$.
As above we get $$\delta(p+aq)+g(b+ta+\langle p,q\rangle)=1.$$ 
Since $a\equiv 1$ modulo $I$, $a=1-x$ for some $x\in I$. 
Let $\beta=(\delta\Phi^{-1})^{*}x$ and $f=xg$. Then it follows that 
$\beta^{*}\Phi=x\delta$. Let $x\delta=\xi$ and  
$\Delta_1=\theta_{(-t,-q)}\tau_{(f,\xi)}\theta_{(t,q)}$. Then 
$\Delta_1 \in {\tran}_{\Sp}((P\perp A^2),I)$. 
Now $\Delta_1 (p,b,a)^t=(p',b',1)^t$ for some $p'\in P$ and $b'\in A$ 
such that $p'\equiv 0$ modulo $I$ and $b'\equiv 0$ modulo $I$. 
So it follows that $\theta_{(p',b')}\equiv \tn{Id}$ modulo $I$. Let 
$\Delta=\theta_{(p',b')}\Delta_1$. 
Then $\Delta \in {\tran}_{\Sp}((P\perp A^2),I)$ and 
$\Delta(p,b,a)^t=(0,0,1)^t$; as required. \hb 
\begin{lm} \label{uni1} 
Let $A$ be an associative ring with identity and 
$P$ and $Q$ be as in \ref{note}.
Let $\Delta$ be a matrix in ${\G}(Q)$. 
If for $m\ge 2$, $\Delta e_m=e_m$, then 
$\Delta\in {\T}(Q) {\G}(P)$. 
\end{lm} 
{\bf Proof.} If $\Delta e_m=e_m$, then in the linear case 
$\Delta$ is of the form 
$$\begin{pmatrix} \beta & 0\\ \gamma & 1\end{pmatrix}
=\begin{pmatrix} 1 & 0 \\ \gamma\beta^{-1} & 1 \end{pmatrix}
\begin{pmatrix} \beta & 0 \\ 0 & 1 \end{pmatrix},$$
for some $\beta\in {\aut}(P)$,  $\gamma\in P^{*}$, and in the 
symplectic case $\Delta$ is of the form
$$\begin{pmatrix} \beta' & p' & 0\\ 0 & 1 & 0 \\
\gamma' & a' & 1\end{pmatrix}
=\begin{pmatrix} 1 & p' & 0 \\ 0 & 1 & 0 \\ 
\gamma'\beta'^{-1} & a' & 1 \end{pmatrix}
\begin{pmatrix} \beta' & 0 & 0 \\ 0 & 1 & 0 \\ 0 & 0 & 1 \end{pmatrix},$$
for some $\beta'\in {\Sp}(P)$, $a'\in A$, $p'\in P$ and 
$\beta_1:A\ra P$ such that $\beta_1(1)=p'$, and 
$\gamma':P\ra A$ chosen in a way that 
$\gamma'=\beta_1^{*}\Phi\beta'$.

Clearly, $\left(\begin{array}{cc} 1 & 0 \\ \gamma\beta^{-1} & 1 
\end{array}\right)$ and 
{\small $\left(\begin{array}{ccc} 1 & p' & 0 \\ 0 & 1 & 0 \\ 
\gamma'\beta'^{-1} & a' & 1 \end{array}\right)=\tau_{(f,\beta_1)}$} are in 
$\tn{ET}(Q)\subset {\T}(Q)$, 
where $f:A\ra A$ given by $1\mapsto a'$. Hence the result follows. \hb

\begin{de} \tn{ Let $k$ be a field. A ring $A$ is said to be essentially of finite type 
over $k$ if $A=S^{-1}C$, with $S$ is a multiplicative closed subset of $C$, and 
$C=k[X_1,\ldots, X_m]/I$ is a quotient ring of a polynomial ring over $k$.}
\end{de}

In (\cite{VOR}, Theorem 3.3), T. Vorst, following ideas of H. Lindel in 
\cite{LH}, proved the following in the linear case:

{\it Let $A$ be a regular local ring essentially of finite type over a perfect 
field $k$. Then  $${\es}_r(A[X])= {\E}_r(A[X])$$ 
for $r\ge 3$.} 

This method of proof also proves the result for the 
symplectic groups. We revisit this proof below. We treat the linear and the 
symplectic cases uniformly. 

\begin{tr} \label{vor}
Let $A$ be a regular local ring essentially of finite type over a perfect 
field $k$. Then  
$${\es}(r,A[X])= {\E}(r,A[X]),$$
for $r\ge 3$ in the linear case, and $r\ge 6$ in the symplectic case.
\end{tr}

We sketch a proof of this theorem below. To prove the theorem we need to use 
the ideas (of A. Suslin and H. Lindel) to establish the statements below in 
the linear case.
 
\begin{lm} \label{vor4} Let $A$  be a commutative ring with identity and 
$S\subset A$ be a multiplicative closed set. If 
${\G}(r,A[X])={\G}(r,A){\E}(r,A[X])$, then 
$${\G}(r,S^{-1}A[X])={\G}(r,S^{-1}A){\E}(r,S^{-1}A[X]).$$ 
\end{lm} 
{\bf Proof.} Let $\alpha(X)\in {\G}(r,S^{-1}A[X])$. Replacing $\alpha(X)$ by 
$\alpha(X)\alpha(0)^{-1}$, we may assume that $\alpha(0)=\tn{Id}$. Let 
$f(X)=\det(\alpha(X)^{-1})$. Then $f(0)=1$. Therefore, 
there exists $s_1\in S$ such 
that $\alpha(s_1X)$ and $f(s_1X)$ are both defined over $A[X]$. Let 
$\alpha_1(X)\in {\G}(r,A[X])$ and $f_1(X)\in A[X]$ with $\alpha_1(0)=I_r$ and 
$f_1(0)=1$, localizing into $\alpha(s_1X)$ and $f(s_1X)$ respectively. Also 
$\det(\alpha(s_1X)).f(s_1X)=1$. Thus, there exists $s_2\in S$ such that 
$\det(\alpha_1(s_2X)).f_1(s_2X)=1$. Hence it follows that 
$\alpha_1(s_2X)\in {\G}(r,A[X])$. Therefore,  
$\alpha_1(s_2X)=\gamma\underset{k=1}{\overset{m}\Pi} ge_{i_k j_k}(f_k(X))$
with $\gamma\in {\G}(r,A)$ and $f_k(X)\in A[X]$ for all $1\le k\le m$. So, 
$$\alpha(X)=\gamma_S \underset{k=1}{\overset{m}\Pi} ge_{i_k j_k}
(f_k (X/\delta_1\delta_2))_S.$$ \hb

We shall assume that $r\ge 3$, in the linear case, and that $r$ is even and 
$r\ge 6$ in the symplectic case.
\begin{pr} \label{vor1} \tn{(A. Suslin)} 
Let $A$  be a commutative ring with identity and $h\in A$ be a
non-nilpotent. Let $\delta\in {\G}(r,A_h)$ and $\sigma(X)=\delta
ge_{kl}(X.f){\delta}^{-1}$, where $k\ne 1$ and $f\in A_h[X]$. Then there
exists a natural number $m$ and a matrix $\tau\in {\E}(r,A[X],XA[X])$
such that $\tau_h=\sigma(h^mX)$.
\end{pr}
{\bf Proof.} For the linear case see (\cite{SUS}, Lemma 3.3). For the 
symplectic case it has been asserted in (\cite{KOP}, $\S$ 3) that a similar 
proof works as in the orthogonal case; and for the orthogonal case 
see (\cite{SUSK}, Lemma 4.6). \hb

\begin{tr} \label{vor3} \tn{(H. Lindel)} 
\tn{(\cite{LH}, Proposition 2 and 3)}
Let $A$ be a regular local ring
essentially of finite type over $k$ with $\dim A\ge 1$, where $k$ is perfect. 
Then there exists a subring $B$ of $A$ with a non-zero divisor $h\in B$ such that
\begin{enumerate}
\item $B$ is the localization of a polynomial ring over $k$,
\item $Ah+B=A$ and $Ah\cap B=Bh$.
\end{enumerate}
\end{tr}
The following was proved by T. Vorst in the linear case in 
(\cite{VOR}, Lemma 2.4 ): 
\begin{lm} \label{vor2} 
Let $A$  be a commutative ring with identity, $B\subset A$, and 
$h\in B$ be a non-nilpotent.
\begin{enumerate}
\item If $Ah+B=A$, then for every $\alpha\in {\E}(r,A_h)$ there exist
 $\beta \in {\E}(r,B_h)$ and $\gamma\in {\E}(r,A)$ such that
$\alpha=\gamma_h\beta$.
\item If moreover $Ah\cap B=Bh$ and $h$ is a non-zero-divisor in $A$, then
for every $\alpha\in {\G}(r,A)$ with $\alpha_h\in {\E}(r,A_h)$ 
there exist a $\beta\in {\G}(r,B)$ and
$\gamma\in {\E}(r,A)$ such that $\alpha=\gamma\beta$.
\end{enumerate}
\end{lm}
{\bf Proof.} $(1)$: Assume that $\alpha=\underset{k=1}{\overset{m}\Pi} 
ge_{i_kj_k}(c_k)$ with $c_k\in A_h$. 
From hypothesis it follows that $Ah^n+B=A$ for all $n$. Hence for all 
$1\le k\le m$ we can find $a_k\in A$, $b_k\in B$ and a natural number $m_k$ 
such that $$c_k=\frac{b_k}{h^{m_k}}+a_kh^s.$$

Let $\sigma_p=
\underset{k=1}{\overset{p}\Pi}ge_{i_kj_k}(c_k)$, $(1\le p\le m)$.  
By Proposition \ref{vor1} there exists a natural number $s$ and 
$\tau_p(X)\in {\E}(r, A[X],XA[X])$ such that 
$$\tau_p(X)=\sigma_p ge_{i_pj_p}(h^sX)\sigma_p^{-1}.$$

So we have 
{\small $$\alpha=\underset{k=1}{\overset{m}\Pi} 
ge_{i_kj_k}\left(\frac{b_k}{h^{m_k}}\right)ge_{i_kj_k}(a_kh^s)= 
\underset{k=1}{\overset{m}\Pi}\sigma_k ge_{i_kj_k}(a_kh^s)\sigma_k^{-1}
\underset{k=1}{\overset{m}\Pi}ge_{i_kj_k}\left(\frac{b_k}{h^{m_k}}\right).$$}
Now let $\gamma=\underset{k=1}{\overset{m}\Pi} \tau_k(a_k)\in {\E}(r,A)$ 
and $\beta= \underset{k=1}{\overset{m}\Pi}ge_{i_kj_k}(\frac{b_k}{h^{m_k}})$. 
Then we are done. \vp 

$(2)$: By hypothesis it follows that $Ah^n\cap B = Bh^n$ for all $n$. 
Hence $B_h\cap A=B$. Using $(1)$ we can write $\alpha_h=\gamma_h\beta$ 
with $\gamma\in {\E}(r,A)$ and $\beta\in {\E}(r,B_h)$. Now 
$\gamma^{-1}\alpha\in
{\G}(r,B)$ and $\beta \in {\G}(r,B_h)$. Moreover 
$(\gamma^{-1}\alpha)_h=\beta$. But this implies that $\gamma^{-1}\alpha 
\in {\G}(r,B)$. Hence $\alpha=\gamma (\gamma^{-1}\alpha)\in 
{\E}(r,A){\G}(r, B).$ \hb \vp \\
{\bf Proof of the Theorem \ref{vor}} We prove the theorem by induction on 
$\dim A$. If $\dim A=0$ then $A$ is a field and the result follows. So we 
assume that $\dim A\ge 1$. 

Let $\alpha(X)\in {\es}(r,A[X])$. 
As the hypothesis of Lemma \ref{vor3} is satisfied, we can find a ring $B$ 
and can choose $h\in B$ as in Lemma \ref{vor3}. Since $\dim A_h < \dim A$,
by induction hypothesis we have that $\alpha_h (X)\in  {\E}(r,A_h[X])$. 
Since $A$ is a regular local ring, 
we have that $h$ is a non-zero-divisor in $A[X]$. 
Now by applying Lemma \ref{vor2} to $\alpha(X)$, we get 
$\alpha(X)=\gamma(X)\beta(X)$ with $\beta(X)\in {\G}(r,A[X])$ and 
$\gamma(X)\in {\E}(r,A[X])$. Hence we have 
$$\alpha(X)=\gamma(X)\gamma(0)^{-1} \beta(0)^{-1}\beta(X),$$
where the first two factors are contained in ${\E}(r,A[X])$. Since the theorem 
is true for a polynomial ring over a field 
(proved in (\cite{SUS}, Corollary 6.7) by A. Suslin for the linear case (and
similarly other cases are also true due to monic inversion) 
and $B$ is a localization of a 
polynomial ring the theorem is also true for $B$ by Lemma \ref{vor4}. Hence
$\beta(0)^{-1}\beta(X)\in {\E}(r,B[X])\subset {\E}(r,A[X])$. \hb 
\begin{tr} \label{pv} 
Let $A$ be a geometrically regular local ring containing a field $k$. Then 
${\es}(r,A[X]) = {\E}(r,A[X])$, for $r\ge 3$, in the linear case, and 
$r\ge 6$, in the symplectic case.
\end{tr} 
{\bf Proof.} If $\dim (A)=0$, then $A$ is a field, and the result 
follows.
Therefore, we assume that 
$\dim (A)\ge 1$. 
In \cite{PD2}, D. Popescu showed that 
if $A$ is a geometrically regular local ring, or when
the characteristic of the residue field is a regular parameter in
$R$, then it is a filtered inductive limit of regular local rings
essentially of finite type over the integers. Hence by Theorem \ref{vor} 
it follows that $${\es}(r,A[X])={\E}(r,A[X])$$ for 
all $r\ge 3$ in the linear case and $r\ge 6$ in the symplectic case. \hb 

\begin{re} \tn{Theorem \ref{pv} is not true for the orthogonal group.
It is not true that ${\es}(r,A) = {\E}(r,A)$, for $r\ge 4$, for 
the orthogonal group, in general, 
even in the case when $A$ is a field. This is known classically due to 
results of Dieudonne, since the spinor norm is surjective.
In the case when $A$ is a local ring similar results have been obtained 
by W. Klingenberg (see \cite{K1}, \cite{K2}), and the 
references therein for the field case.} 
\end{re}

\begin{re} \tn{The proof of Theorem \ref{pv} can be used to show that if 
 $A$ is a geometrically regular local ring containing a field $k$ then 
a stably elementary orthogonal matrix $\sigma(X) \in {\rm SO}_{2n}(A[X])$, 
$n \geq 3$, with $\sigma(0) = I_{2n}$, is an elementary orthogonal matrix.}
\end{re}

\begin{re} \tn{Using ``deep splitting technique'' as in 
(\cite{Blr}, Definition 3.6, Corollary 3.9) 
one can show that Lemma \ref{vor2} is valid for $r=4$. Consequently, 
Theorem \ref{vor} and Theorem \ref{pv}, are also valid for $r=4$. The 
above remark is also true when $n = 2$.} 
\end{re} \vp 

\medskip

We now establish the main theorems stated in the Introduction.

~~~To have a uniform notation in Theorem \ref{st1} we use the notation 
$\widetilde{{\es}}(P)$ to denote ${\SL}(P)$ and ${\Sp}(P)$ and 
$\widetilde{{\T}}(P)$ to denote ${\tran}(P)$ and ${\tran}_{\Sp}(P)$.  
\vp \\
{\bf Proof of Theorem 1 and 2.} 
The homotopy technique used here is as in (\cite{RV}, Proposition 3.4).
In view of L.N. Vaserstein's result in \cite{V2}, to prove the result it is 
enough to prove the injectivity for $n=d+1$. 
Let $n_1=n+1$ and $n+2$ in the linear and the symplectic cases respectively.
Let $\widetilde{Q}$ denote $P\oplus A$ in the linear case and $P\perp A^2$ in 
the symplectic case. Consider $\gamma\in \widetilde{{\es}}(P)$  such that 
$\widetilde{\gamma}=\gamma\perp \tn{Id}
\in \widetilde{{\T}} (Q)$.  
Let $\eta(X)$ be the isotopy between $\widetilde{\gamma}$ and identity. As
before, viewing $\eta(X)$ as a matrix (as in \ref{endo}), it follows that 
$v(X)^t$, where $v(X)=\eta(X)e_{n_1}$, is a unimodular element in 
$\widetilde{Q}[X]$.
Note that $v(X)\equiv e_{n_1}$ modulo $(X^2-X)$. 
Using Lemma \ref{rv} for Theorem 1 and Lemma \ref{b3a} for Theorem 2
over $A[X]$ it follows that there exists
$\sigma(X)\in \widetilde{{\es}}(\widetilde{Q}[X])$ 
such that $\sigma(X)^tv(X)=e_{n_1}$ and $\sigma(X)\equiv $ Id modulo 
$(X^2-X)$. Therefore, $\sigma(X)^t\eta(X)e_{n_1}=e_{n_1}$. 
Hence by Lemma \ref{uni1}, 
$\sigma(X)^t\eta(X)=\xi(X)\widetilde{\eta}(X),$ where 
$\xi(X)\in \widetilde{{\E}}(\widetilde{Q}[X])$  
and $\widetilde{\eta}(X)\in \widetilde{\es}(P[X])$. 
Since $\sigma(X)\equiv $ Id modulo 
$(X^2-X)$, $\widetilde{\eta}(X)$ is an isotopy between $\gamma$ and 
the Identity. 

Now assume $A$ is regular and contains a field $k$. Hence 
for every prime ideal $\mf{p}\in {\sp}(A)$, 
$$\widetilde{\eta_{\mf{p}}}(X)\in \begin{cases} \widetilde{\es}(n,A_{\mf{p}}[X])=
\widetilde{\E}(n,A_{\mf{p}}[X]) \mbox{ in the linear case, and }\\
\widetilde{\es}(2n,A_{\mf{p}}[X])=
\widetilde{\E}(2n,A_{\mf{p}}[X]) \mbox{ in the symplectic case.}
\end{cases} $$ 
(by Theorem \ref{pv}). 
Since $\widetilde{\eta}(0)$=Id, by the L-G Principle (Theorem \ref{lgt}) 
for the tranvection groups  we get 
$\widetilde{\eta}(X)\in \widetilde{{\ET}}(P[X])$.
 Whence $\gamma=\widetilde{\eta}(1)\in \widetilde{{\ET}}(P)$; as required. 
$\tn{ }$ \hb

\noindent
{\bf Acknowledgement:} The authors thank Professor R.G. Swan profusely for 
his many illuminating comments and corrections, and for his unstinting 
support and encouragement. 


\addcontentsline{toc}{chapter}{Bibliography} \vp

\begin{thebibliography}{99} 
\setlength{\itemsep}{-.5ex}
\bibitem{BR2} A. Bak, R. Basu, Ravi.A. Rao; 
Local-Global Principle for Transvection Groups, preprint, see  
{\it  http://arxiv.org/abs/0908.3094}

\bibitem{brj} R. Basu, Ravi A. Rao, Selby Jose: Injective Stability for ${\k}$ 
of the Orthogonal group, provisionally accepted in Journal of Algebra. 
\bibitem{B3} H. Bass, J. Milnor, J. -P. Serre; Solution of the congruence 
subgroup problem for ${\SL}_n$, $n\ge 3$ and ${\Sp}_{2n}$, $n\ge 2$.  
Inst. Publ. Math. IHES 33 (1967), 59--137.
\bibitem{B2} H. Bass; Algebraic K-theory, Math. Lecture note series, 
W.A. Benjamin, Inc. (1968). 
\bibitem{Bh} S.M. Bhatwadekar; A cancellation theorem for projective 
modules over affine algebras over $C_1$-fields. Journal of Pure and 
Applied Algebra 183 (2003), no. 1-3, 17--26. 
\bibitem{Blr} S.M. Bhatwadekar, H. Lindel, R.A. Rao; 
The Bass-Murthy question: Serre dimension of Laurent polynomial extensions.  
Invent. Math.  81  (1985),  No. 1, 189--203.
 \bibitem{K2} Wilhelm Klingenberg; Orthogonale Gruppen �ber lokalen Ringen. 
(German) American J. Math.  83  (1961), 281--320. 
\bibitem{K1} Wilhelm Klingenberg; Orthogonal groups over local rings. 
Bull. Amer. Math. Soc. 67 (1961), 291--297. 
\bibitem{KOP}  V.I. Kope\u\i ko; The stabilization of symplectic groups 
over a polynomial ring. Math. USSR. Sbornik 34 (1978), 655--669.
\bibitem{LAM} T. Y. Lam; Serre's Conjecture. Lecture Notes in
Mathematics, Vl. 635. Springer-Verlag, Berlin-New York, 1978. 
\bibitem{LH} H. Lindel; On the Bass-Quillen conjecture concerning
projective modules over polynomial rings. Invent. Math.
65 (1981--82), no. 2, 319--323. 
\bibitem{PD2} D. Popescu; Polynomial rings and their projective
modules, Nagoya Math. J. 113 (1989), 121--128. 

\bibitem{RV} Ravi.A. Rao, Wilberd van der Kallen; Improved stability for 
${\k}$ and $\tn{WMS}_ d$ of a non-singular affine algebra. ${\K}$-theory 
(Strasbourg, 1992).  Asterisque  no. 226 (1994), 11, 411--420.

\bibitem{su1} A.A. Suslin; Stably Free Modules. (Russian) Math. USSR Sbornik
102 (144) (1977), no. 4, 537--550.
Mat. Inst. Steklov. (LOMI)  114  (1982), 187--195, 222.
\bibitem{SUS} A.A. Suslin; On the structure of special linear group over
polynomial rings. Math. USSR. Izv. 11 (1977), 221--238.
\bibitem{SUSK}  A.A. Suslin, V.I. Kopeiko; Quadratic modules and orthogonal 
groups over polynomial rings. Nauchn. Sem., LOMI 71 (1978), 216--250.
\bibitem{su2} A.A. Suslin; Cancellation for affine varieties. (Russian) Modules and algebraic groups.  Zap. Nauchn. Sem. Leningrad. Otdel. Mat. Inst. Steklov. (LOMI)  114  (1982), 187--195, 222.
\bibitem{S} R.G. Swan, Serre's Problem. 
Conference on Commutative Algebra--1975, 
Queen's Paper in Pure and Applied Mathematics, no. 42, (1975), 1--60. 
\bibitem{V} L.N. Vaserstein; On the stabilization of the general linear group 
over a ring.  Mat. Sbornik (N.S.) 79 (121) 405--424 (Russian); English 
translated in Math. USSR-Sbornik. 8 (1969), 383--400. 
\bibitem{V2} L.N. Vaserstein; Stabilization of Unitary and Orthogonal Groups 
over a Ring with Involution. Mat. Sbornik, Tom 81 (123) (1970) no. 3, 307--326.
\bibitem{V1} L.N. Vaserstein; Stable range of rings and dimensionality of 
topological spaces. Fuct. Anal. Appl. 5 (1971), no. 2, 102--110.
\bibitem{V3} L.N. Vaserstein; Stabilization For Classical Groups over Rings. 
(Russian)  Mat. Sb. (N.S.)  93 (135)  (1974), 268--295, 327. 
\bibitem{VOR}  T. Vorst; The general linear group of polynomial rings over 
regular rings. Comm. Algebra 9 (1981) no. 5, 499--509.
\end{thebibliography}
\end{document}